\documentclass[a4paper,12pt]{amsart}

\usepackage{amsfonts,latexsym,rawfonts,amsmath,amssymb,amsthm, mathrsfs, lscape}
\usepackage{newlfont} 
\usepackage[all]{xy}
\usepackage[english]{babel}
\usepackage[utf8]{inputenc}
\usepackage[margin=3.5cm]{geometry}
\usepackage{color, hyperref, comment}
\usepackage{nicematrix}

\newtheorem{lemma}{Lemma}[section]
\newtheorem{theorem}[lemma]{Theorem}
\newtheorem{proposition}[lemma]{Proposition}
\newtheorem{conjecture}[lemma]{Conjecture}

\newtheorem{remark}[lemma]{Remark}
\newtheorem*{claim*}{Claim}
\newtheorem{nota}[lemma]{Notation}

\title{On the dimension of affine subspaces of nilpotent matrices}

\author{Simone Calamai}
\address[Simone Calamai]{Dipartimento di Matematica e Informatica ``Ulisse Dini''\\
Università di Firenze\\
viale Morgagni 67/A\\
50134 Firenze, Italia}
\email{simone.calamai@unifi.it}

\author{Elena Rubei}
\address[Elena Rubei]{Dipartimento di Matematica e Informatica ``Ulisse Dini''\\
Università di Firenze\\
viale Morgagni 67/A\\
50134 Firenze, Italia}
\email{elena.rubei@unifi.it}

\keywords{affine subspaces, nilpotent matrices}
\subjclass[2020]{15A30}

\newcommand{\N}{\mathbb{N}}

\begin{document}

\begin{abstract}
The focus of the paper is on the maximal dimension of affine subspaces 
of nilpotent $n \times n $ matrices with fixed rank.
In particular we obtain a result in the "border" case rank equal to $n-1$; we also consider the much simpler case rank equal to $1$.
\end{abstract}

\maketitle

\section{Introduction}

There is a wide literature on the maximal dimension of linear or affine subspaces
of matrices with specific characteristics.
The most famous result on linear subspaces of nilpotent matrices is the following theorem.

\begin{theorem} \label{Gerst} [Gerstenhaber, Serezhkin] Let $K$ be  a field. The maximal dimension of a linear subspace of nilpotent $n \times n$  matrices over $K$  
 is $\frac{n(n-1)}{2}$. 
\end{theorem}

The result above  was proved by Gerstenhaber  under the assumption that $K$ has at least $n$ elements (see \cite{Gerst}) and in 1985  the theorem was generalized by Serezhkin for any field (see \cite{Ser}). Moreover, Gerstenhaber proved that a subspace attains the maximal dimension if and only if it is similar to the space of all strictly upper-triangular matrices.
We also  mention  that, in 1991, Mathes, Omladi\u{c} and Radjavi 
gave  new simple proofs of Theorem \ref{Gerst} and of the result on the subspaces attaining the maximal dimension under the condition that the field has more than two elements (see  \cite{M-O-R}, Theorems 1 and 2).

In 1959, Gerstenhaber proved that, if $K$ is a field with  at least $r+1$ elements, then the dimension of a linear subspace of nilpotent $n \times n$ matrices over $K$  with rank at most $r$ is less than or equal to  $$\frac{n(n-1)}{2} - \frac{(n-r)(n-r-1)}{2} $$ (see \cite{Gerst1}, Theorem 2.7) and in 1962 he gave an improved bound on the dimension of a linear subspace of nilpotent matrices in terms of the sizes of the Jordan blocks in the Jordan forms of the elements of the subspace. To state Gerstenhaber's result we need to fix some notation: for a positive integer $n$, we say that $(a_1,\ldots, a_n)$ is a partition of $n$ if $a_1 \geq a_2 \geq \ldots \geq a_n \geq 0$ and $n=a_1+ \ldots + a_n$; the conjugate partition of $(a_1,\ldots, a_n)$  is the partition  $(\bar{a}_1,\ldots, \bar{a}_n)$ of $n$, where $\bar{a}_j$ is the number of the $a_i$ that are greater than or equal to $j$; given two partitions of $n$,
$(a_1,\ldots, a_n)$ and $(b_1,\ldots, b_n)$, we say that $(a_1,\ldots, a_n) \leq (b_1,\ldots, b_n)$ if  $a_1+ \ldots +a_j \leq b_1+ \ldots +b_j$ for any $j=1,\ldots ,n$; finally,
for any nilpotent $n \times n $ matrix $A$, let $k_1, \ldots, k_l$ with $k_1 \geq k_2 \geq \ldots \geq k_l \geq 1$ be the sizes of the Jordan blocks in the Jordan form of $A$; the partition $(k_1, \ldots , k_l ,0 \ldots,0 )$, where $0$ is repeated $n-l$ times, is called the Jordan partition of $n$ associated to $A$. The result of Gerstenhaber linking the dimension of a linear subspace of nilpotent matrices to the sizes of the Jordan blocks is the following, see \cite{Gerst2}: 

\begin{theorem} \label{Gerst2} [Gerstenhaber]  Let $K$ be a field. Let ${\cal S}$ be a linear subspace of nilpotent $ n \times n $  matrices. Let $(a_1, \ldots, a_n)$ be the least upper bound of the set of all Jordan partitions of matrices in $S$. If $|K|$ is sufficiently large, then 
$$ \dim({\cal S}) \leq \frac{1}{2} \left(  n^2- \sum_{i=1}^{n} \bar{a}_i^2 \right).$$
\end{theorem}

We also point out  that recently Quinlan and De Seguins Pazzis generalized Theorem \ref{Gerst} proving that the maximal dimension of a linear subspace of $n \times n$ matrices  over a field $K$ with no nonzero eigenvalue in $K$ is $\frac{n(n-1)}{2}$ (see  \cite{Q}, Theorem 1.3 and \cite{DP2}, Theorem 9).


In 1993,  Brualdi and Chavey gave another proof of Theorem \ref{Gerst2} (see Theorem 3.1 in  \cite{B-C}) and gave a bound on the maximal dimension of a linear space of nilpotent matrices with bounded nilindex when the cardinality of the field is sufficiently large, where for the nilindex of a nilpotent matrix $A$ we mean the minimum $k$ such that $A^k$ is zero (see  Corollary 3.6 in  \cite{B-C}).

More recently there has been a significant progress  on this area thanks to the papers of MacDonald, MacDougall and Sweet.

In 2009, MacDougall and Sweet proved that, if $|K| >n/2$, the maximal dimension of  a linear space of nilpotent $n \times n$ matrices over $K$ with nilindex less than or equal to $2$ and rank less than or equal to $r$ is $r(n-r)$ (see  \cite{S-M}, Theorem 4 and \S4).

In 2012, MacDonald, MacDougall and Sweet proved the following result (see \cite{M-M-S}, Theorems 5.1 and 6.1):

 \begin{theorem} \label{MMS-theorem}[MacDonald, MacDougall, Sweet]
 If $K$ is a field with $|K| >n$ and ${\cal S}$ is a  linear subspace of nilpotent  $n \times n$ matrices whose nilindex is less than or equal to $k$ and rank less than or equal to $r$, then $$ \dim({\cal S}) \leq nr - \frac{r^2}{2} - \frac{r}{2}  + \frac{q^2}{2} (k-1) +
 \frac{q}{2} (-2r +k-1),$$
where $q=  \left\lfloor \frac{r}{k-1} \right\rfloor$. The bound above
is sharp.
\end{theorem}

Moreover, they improved Gerstenhaber's bound on the 
  dimension of a linear subspace of nilpotent matrices in terms of the sizes of the Jordan blocks in the Jordan forms of the elements of the subspace, see Theorem 7.2 in \cite{M-M-S} for the precise statement.

Finally,  we also mention  that   Omladi\u{c} proved (see \cite{O}, Corollary 8) that for  a maximal linear subspace of nilpotent $n \times n$ matrices  the following conditions are equivalent: (a) it has dimension $\frac{n(n-1)}{2}$, (b) it is triangularizable, (c)  it contains a matrix of rank $n-1$ and its square.

In this paper, we focus on the maximal dimension of affine subspaces 
of nilpotent $n \times n$ matrices with fixed rank.

We mention that in \cite{Rubei2024specific} (see Theorem 3) it has been proved that
the maximal dimension of an affine subspace of nilpotent $n \times n$ matrices over any field $K$
is $ \frac{n(n-1)}{2}$ and that, if the characteristic of $K$ is zero, 
an affine not linear  subspace of nilpotent $n \times n$ matrices has dimension  less than or equal to $ \frac{n(n-1)}{2}-1$. 

Here we obtain a result on  the maximal dimension of affine subspaces of nilpotent $n \times n$ matrices with 
 rank equal to $n-1$ (and this condition is obviously equivalent to the condition that the nilindex is equal to $n$):

\begin{theorem} \label{dimker=1}
Let $n \in \N-\{0\}$ and
let $K$ be a field such that $|K| \geq n+1$.
Then the maximal dimension of an affine subspace of nilpotent  $n \times n$ matrices over $K$ with rank equal to $n-1$  is $$ \frac{(n-1)(n-2)}{2}.$$
\end{theorem}

Observe that the result  we have found for the maximal dimension of an affine space of nilpotent $n \times n$ matrices of constant rank $n -1$ is substantially different from the case of linear spaces; in fact, given a linear space of  nilpotent $n \times n$ matrices of constant rank $n -1$, we have that its dimension   can be at most $2$ and $n$ must be odd, as proved by Causa, Re and Teodorescu in  \cite{C-R-T}, Theorem 2.3
(we recall that  a linear subspace of matrices is of constant rank $r$ if every of its nonzero elements   has rank $r$).

Moreover, in Proposition \ref{rank=1}, we consider the other border case, that is, the case rank equal to $1$; 
this case is much simpler than the case rank equal to $n-1$ and
it is easy to prove 
 that, if the cardinality of the field $K$ is at least $3$,  the maximal dimension of an affine subspace of nilpotent  $n \times n$ matrices over $K$ with rank equal to $1$  is $n-2$.

At the end of the paper we have included a conjecture  on the maximal dimension of affine subspace of nilpotent matrices of fixed  rank that could generalize the results above (see Conjecture \ref{conjecture}).

Finally, we point out that, obviously,  the problem could have some links  with the study of the maximal dimension of linear or affine subspaces of matrices with constant rank or bounded rank,
 which has attracted wide literature, see for instance  \cite{B-F-M}, \cite{DP3}, \cite{DP4},
  \cite{I-L},  \cite{L-M}, \cite{M-Mez}, \cite{M-M}, 
\cite{Rubei2022},  \cite{Rubei2024antisymm},   \cite{W1} for papers on 
the maximal dimension of linear or affine subspaces of matrices with constant rank and the 
papers \cite{A}, \cite{DP2}, \cite{DP3}, \cite{E-H}, \cite{Fl}, \cite{H-L},  \cite{Rubei2024range} for the problem
of the maximal dimension of a linear or an affine subspace of matrices whose rank is bounded below or above or both below and above and the references therein.

\section{Proofs of the results}

\begin{nota} \label{notaz} 

By natural numbers we mean non-negative integers and 
we denote the set of  natural numbers  by $ \mathbb{N}$.

Let $m,n \in \mathbb{N} -\{0\} $ and $K$ be a field.

We denote the $n \times n$ identity matrix over $K $ by $I^K_n$ and
the $m \times n$ null matrix over $K$ by $0^K_{m \times n }$. We define $J^K_n$ 
to be the $n \times n$  matrix over $K $ such that 
$$ (J^K_n)_{i,j}= \left\{  \begin{array}{ll}  
1 &  \mbox{\rm if} \ \; j=i+1 \\
0 &  \mbox{\rm otherwise .}
\end{array} \right.$$

Moreover, we denote by $E_{i,j}^{n,K}$ the $n \times n$  matrix over $K$ such that 
$$ (E_{i,j}^{n,K})_{x,y} = \left\{ \begin{array}{ll}
1 & \mbox{\rm if} \; (x,y)=(i,j) \\
0 & \mbox{\rm otherwise.}
\end{array} \right.$$
We define $e_i^{n,K}$ to be the $i$-th element of the canonical basis of $K^n$.

 We omit the superscripts and the subscripts in  $I^K_n$, $0^K_{m \times n}$ and  $J^K_n$ and  
  the superscripts in $E_{i,j}^{n,K}$ and  $e_i^{n,K}$ 
 when the underlying field and the size of the matrix (or of the vector in the last case) are clear from the context.

For any $A \in  M(m \times n,K) $, the submatrix of $A$ given by the rows $i_1, , \ldots , i_k$ and the columns $j_1, \ldots, j_s$  is denoted by $A^{(j_1, \ldots, j_s)}_{(i_1, \ldots, i_k)}$.

For any matrix $A$, we denote the rank of $A$ by  $rk(A)$  and, for any square matrix $A$,  we denote   the trace of $A$ by $tr(A)$.


For any subset ${\cal S}$ of a vector space ${\cal V}$, the linear subspace generated by the elements of ${\cal S}$ is denoted by $\langle {\cal S} \rangle $.

Let ${\cal V}$ be a vector space; let $a \in {\cal V}$ and ${\cal Z}$ be a vector subspace of  ${\cal V}$; consider the affine subspace  ${\cal S} := a +{\cal Z}$; we call ${\cal Z}$ the direction of ${\cal S}$. 

\end{nota}

In \cite{M-O-R} the following proposition is proved (see Corollary 1) :

\begin{proposition} \label{MOR} [Mathes, Omladi\u{c}, Radjavi]
Let ${\cal V}$ be a linear space of nilpotent matrices over a field $K$ of characteristic $0$. If $A, B \in {\cal V}$ and $m \geq 1$ , then $tr(A^mB)=0$.
\end{proposition}

In 2012 MacDonald, MacDougall and Sweet generalized the proposition above in the following way  (see  \cite{M-M-S}, Theorem 2.4):

\begin{proposition} \label{MMS} [MacDonald, MacDougall, Sweet]  Let $m \in \N$.
Let ${\cal V}$ be a linear space of nilpotent matrices over a field $K$ such that $|K| >m$. If $A, B \in {\cal V}$, then $tr(A^mB)=0$.
\end{proposition}

\begin{lemma} \label{J}
Let $n \in \N-\{0\}$ and $K$ be a field.

(i) Let $A \in M(n \times n, K)$. Then 
$JA$ is the matrix whose rows are $A_{(2)}, \ldots, A_{(n)}, 0 $ and $AJ $ is the matrix whose columns are $0, A^{(1)}, \ldots, A^{(n-1)}$.

(ii) For any $c_1, \ldots, c_{n-1} \in K$, the matrix $I+ \sum_{h=1}^{n-1} c_h J^h$
is invertible and the inverse is a matrix of the kind $I+ \sum_{h=1}^{n-1} d_h J^h$ for some $d_1, \ldots, d_{n-1} \in K$.

(iii) Let $A \in M(n \times n, K)$ and  $C=I+ \sum_{h=1}^{n-1} c_h J^h$ for some $c_1, \ldots, c_{n-1} \in K$. Let $i, j \in \{1, \ldots, n\}$; if $$A_{a,b}=0 \;\; \; \;\forall (a,b) \in \{i, \ldots, n\} \times  \{1, \ldots, j\} -\{(i,j)\},$$
then $$(C^{-1}AC)_{i,j}=A_{i,j}.$$

(iv) Let $A \in M(n \times n, K)$ and $l \in \{2, \ldots, n\}$. If $$A_{l,1} \neq 0 \;\; \mbox{\rm and} \;\; A_{i,1}=0 \;\forall i \in \{l+1, \ldots n\},$$ then there exist $c_1, \ldots, c_{n-1} \in K$ such that, if we define  $C= I+ \sum_{h=1}^{n-1} c_h J^h$, we have that   
$$(CAC^{-1})_{i,1}=0 \;\; \; \;\forall i \in \{1, \ldots, l-1\}.$$

\end{lemma} 

\begin{proof}
The first statement is obvious. As to (ii), it follows immediately from the fact that, if $N$ is a nilpotent $n \times n $ matrix, then $I+N$ is invertible and the inverse is $I-N+N^2 -\ldots +(-1)^{n-1} N^{n-1}$. 

The statement (iii) follows  from (i) and (ii), in fact: if we write  $C^{-1}=I+ \sum_{h=1}^{n-1} d_h J^h$ for some $d_1, \ldots, d_{n-1} \in K$ (see part (ii) of the lemma we are proving), then 
$$(C^{-1}AC)_{i,j}= \left( \left(I+ \sum_{a=1}^{n-1} d_a J^a\right) A \left(I+ \sum_{b=1}^{n-1} c_b J^b\right) \right)_{i,j}= $$ $$=\left(A+ \sum_{a=1}^{ n-1} d_a J^a A  + \sum_{b=1}^{n-1} c_b AJ^b+ 
 \sum_{a=1}^{ n-1} \sum_{b=1}^{n-1}
 d_a c_b J^a A J^b \right)_{i,j}=
A_{i,j},$$
where in the last equality we used part (i) of  the lemma we are proving.

Finally, to prove (iv), observe that 
if we write, as before,  $C^{-1}=I+ \sum_{h=1}^{n-1} d_h J^h$ for some $d_1, \ldots, d_{n-1} \in K$, then 
$$(CAC^{-1})_{i,1}= \left( \left(I+ \sum_{a=1}^{n-1} c_a J^a\right) A \left(I+ \sum_{b=1}^{n-1} d_b J^b\right) \right)_{i,1}= $$ 
$$=\left(A+ \sum_{a=1}^{n-1} c_a J^a A  +
 \sum_{b=1}^{ n-1} d_b AJ^b+   
 \sum_{a=1}^{ n-1} \sum_{b=1}^{n-1}
   c_a d_b J^a A J^b \right)_{i,1}=$$ 
$$=\left(A+ \sum_{a=1}^{ n-1} c_a J^a A \right)_{i,1}=
A_{i,1} +  \sum_{a=1}^{ n-1} c_a  A_{i+a,1} ,$$
hence 
$$ (CAC^{-1})_{l,1}= A_{l,1} +  \sum_{a=1}^{ n-1} c_a  A_{l+a,1}  =  A_{l,1} $$
$$ (CAC^{-1})_{l-1,1}= A_{l-1,1} +  \sum_{a=1}^{ n-1} c_a  A_{l-1+a,1}  =  A_{l-1,1} + c_1 A_{l,1} $$
$$ (CAC^{-1})_{l-2,1}= A_{l-2,1} +  \sum_{a=1}^{n-1} c_a  A_{l-2+a,1}  =  A_{l-2,1} + c_1 A_{l-1,1} + c_2 A_{l,1}  $$
and so on. 
From the formulas above, we can see easily that we can choose the $c_i$ as we want: take $c_1$ such that $(CAC^{-1})_{l-1,1} = 0$, then take $c_2$ such that  $(CAC^{-1})_{l-2,1} = 0$ and so on.
\end{proof}

\begin{lemma} \label{primorem}
Let $n \in \N-\{0\}$ and
let $K$ be a field such that $|K| \geq n+1$. 

Let ${\cal S}$ be an affine subspace of $M(n \times n, K)$ such that every element of ${\cal S}$ is nilpotent and let ${\cal Z}$ be its direction (see Notation \ref{notaz} for the definition of direction). Then also every element of ${\cal Z}$ is nilpotent. 

Moreover, if ${\cal S} = J +{\cal Z}$, where $J$ is defined in Notation \ref{notaz}, we have that  $A_{n,1}=0 $ for every $A \in {\cal Z}$.

\end{lemma}

\begin{proof}
Let ${\cal S}= P+{\cal Z}$ for some $P \in M(n \times n, K)$;
then for every $A \in {\cal Z}$, we have that $(P+tA)^n=0$ for every $t \in K$,
hence, 
for every $i,j \in \{1,\ldots, n\}$ and $t \in K$, $\left((P+tA)^n\right)_{i,j}=0$. Obviously, 
 $\left((P+tA)^n\right)_{i,j}$ is a polynomial in $t$ of degree less than or equal to $n$ whose term of degree $n$ is $t^n (A^n)_{i,j}$;
since $|K| \geq n+1 $ there must exist an element in $K$ that is not a root of this polynomial   unless every of its 
coefficients is zero; hence $(A^n)_{i,j} =0 $  for every $i,j \in \{1,\ldots, n\}$; so $A$ is nilpotent.

Suppose now that  ${\cal S} =J+{\cal Z}$; let $A \in {\cal Z}$; by the nilpotency of the elements of ${\cal S}$, we have 
 $\det (J+t A)=0$ for every $t \in K$;
but $|K| \geq n+1 $ and $\det (J+t A)$ is a polynomial in $t$ of degree at most $n$, so it can be zero for every element of $K$ if and only if each of its coefficients  is zero. Since the coefficient of the term of degree $1$ is $\pm A_{n,1}$, we can conclude. 
\end{proof}

\begin{lemma} \label{triangblocchi}
Let $n \in \N$ with $n \geq 2$ and
let $K$ be a field.
Let ${\cal S}=J+{\cal Z}$ be an affine subspace of $M(n \times n, K)$ such that every element of ${\cal S}$ is nilpotent and of rank $n-1$, where $J$ is defined in Notation \ref{notaz} and ${\cal Z}$ is a vector subspace. 

Suppose that there exists $r \in \{1, \ldots,n-1\} $  such that every element of ${\cal Z}$ is an upper triangular block matrix whose diagonal blocks are of size $r \times r $ and $ (n-r) \times (n-r)$.

If for every  $n' \in \{1, \ldots,n-1\} $ 
the maximal dimension of an affine subspace of nilpotent  $n' \times n'$ matrices over $K$ with rank equal to $n'-1$ 
is less than or equal to 
$ \frac{(n'-1)(n'-2)}{2},$ then $\dim ({\cal Z}) \leq  \frac{(n-1)(n-2)}{2}$.
\end{lemma}

\begin{proof}
Let $M\in {\cal Z}$. By assumption there exist $A\in M(r \times r, K)$, $B\in M(r \times (n-r), K)$, $C\in M((n-r) \times (n-r), K)$ such that 
$ M = \left(
\begin{array}{c|c}
   A & B \\  \hline
 \mathbf{0} &  C 
\end{array} \right),
$ and hence
$$
J_n + M =  \left(
\begin{array}{c|c}
  J_r + A & E_{r,1} + B \\
  \hline
 \mathbf{0} & J_{n-r} + C 
\end{array} \right)
.
$$
We observe that both $J_r + A$ and $J_{n-r} + C$ are nilpotent since $J_n+ M$ is a nilpotent upper triangular block matrix. 

Moreover, we claim that $rk(J_r + A)=r-1$ and $rk(J_{n-r} + C)= n-r-1$. In order to obtain $rk(J_r + A)=r-1$, first observe that the nilpotency of  $J_r+A$  implies  $rk(J_r + A)\leq r-1$. Now, if $r=1$ the only possibility is that $rk(J_r + A) = r-1$ and we are done. 
If $r\geq 2$  and we assume by contradiction that  $rk(J_r + A) < r-1$, 
 the $n\times r$ submatrix of $J_n+M$ given by the first $r$ columns   would have at most $r-2$ independent columns. This  would imply that 
 $J_n +M$ has at most $n-2$ independent columns,  which contradicts the assumption that each element of ${\cal S}=J_n +{\cal Z}$ has rank exactly $n-1$. Thus, we proved the claim that $rk(J_r + A)=r-1$; the proof of the claim $rk(J_{n-r} + C)= n-r-1$ is similar. 

We define the map
$$
 \begin{aligned}
    \alpha : \phantom{AAAAAA} {\cal Z} \phantom{AA} & \longrightarrow  M(r \times r , K) \\
    M= \left(
\begin{array}{c|c}
   A &  B \\
  \hline
 \mathbf{0} & C 
\end{array} \right) & \longmapsto A \; .
  \end{aligned}
$$
We define ${\cal S}_r := J_r + \alpha({\cal Z})$ and we notice that, for what we proved at the beginning of the argument, ${\cal S}_r$ is an affine subspace of $M(r\times r, K)$ such that every of its elements is nilpotent and has rank $r-1$. Thus, 
by assumption (take $n'=r$), we get that $\dim\mbox{Im}(\alpha) \leq \frac{(r-1)(r-2)}{2}$. 

Next, we define the map 
$$
 \begin{aligned}
    \gamma : \phantom{AAAA}\mbox{Ker}(\alpha)\phantom{AA} & \longrightarrow  M((n-r) \times (n-r) , K) \\
    M= \left(
\begin{array}{c|c}
   \mathbf{0} &  B \\
  \hline
 \mathbf{0} & C 
\end{array} \right) & \longmapsto C \;\; . 
  \end{aligned}
$$
Again ${\cal S}_{n-r}:= J_{n-r}+\gamma(\mbox{Ker}(\alpha))$ is an affine subspace of $M((n-r)\times(n-r),K)$ such that every of its elements is nilpotent and of rank $n-r-1$; thus, by assumption (take $n'=n-r$), we get $\dim\mbox{Im}(\gamma) \leq \frac{(n-r-1)(n-r-2)}{2}$. In summary, so far we have:
$$
\begin{aligned}
    \dim ( {\cal Z}) &= \dim\mbox{Im}(\alpha) + \dim\mbox{Ker}(\alpha) \leq
\frac{(r-1)(r-2)}{2}+\dim\mbox{Ker}(\alpha) =\\
&=\frac{(r-1)(r-2)}{2}+\dim\mbox{Im}(\gamma) +\dim\mbox{Ker}(\gamma) \leq \\
&\leq \frac{(r-1)(r-2)}{2} + \frac{(n-r-1)(n-r-2)}{2} +\dim\mbox{Ker}(\gamma) \; .
\end{aligned}
$$
Regarding $\mbox{Ker}(\gamma)$, we observe that if $M\in \mbox{Ker}(\gamma)$ then $M= \left(
\begin{array}{c|c}
   \mathbf{0} &  B \\
  \hline
 \mathbf{0} & \mathbf{0} 
\end{array} \right)$, and we claim that it must hold that $B_{r,1} = 0$. In fact, if  by contradiction $B_{r,1} \neq 0$, then there would exist $s\in K$ such that $E_{r,1}+sB_{r,1}=0$, and for such $s$ it would hold that 
$$
J_n + sM = \left(
\begin{array}{c|c}
   J_r &  E_{r,1}+sB \\
  \hline
 \mathbf{0} & J_{n-r}
\end{array} \right)
=\left(
\begin{array}{c|c}
   \mathbf{0} &  W \\
  \hline
 \mathbf{0} & \mathbf{0} \, 
\end{array} \right),
$$
where in the last equality we modified the sizes of the blocks in such a way that $W\in M((n-1)\times (n-1), K)$. About the matrix $W$ we notice that it is upper triangular, and that its diagonal is
$$
(
\underbrace{1,\,\ldots,\,1}_\text{(r-1)-times} ,\, 0 ,\, \underbrace{1,\,\ldots,\,1}_\text{(n-r-1)-times})
$$
whence $rk(W)= n-2$, which would imply 
$rk(J_n + sM) = n-2$ that contradicts 
the fact that each element of $J_n + {\cal Z}$ has rank $n-1$. Thus we proved the claim that, if $M\in \mbox{Ker}(\gamma)$, then $M= \left(
\begin{array}{c|c}
   \mathbf{0} &  B \\
  \hline
 \mathbf{0} & \mathbf{0} 
\end{array} \right)$ with $B_{r,1} = 0$. This tells us that $\dim\mbox{Ker}(\gamma)\leq r(n-r)-1$; with such inequality we conclude the wanted estimate of $\dim ({\cal Z})$ as follows
:$$
\begin{aligned}
    \dim ({\cal Z}) &\leq \frac{(r-1)(r-2)}{2} + \frac{(n-r-1)(n-r-2)}{2} +\dim\mbox{Ker}(\gamma) \leq \\
    &\leq  \frac{(r-1)(r-2)}{2} + \frac{(n-r-1)(n-r-2)}{2} + r(n-r)-1 =\\
    &= \frac{(n-1)(n-2)}{2} \; .
\end{aligned}
$$
\end{proof}

\begin{lemma} \label{dimL} Let $n \in \N$ with $n \geq 3$ and
let $K$ be a field such that $|K| \geq n+1$. 
Let ${\cal S}=J+{\cal Z}$ be an affine subspace of $M(n \times n, K)$ such that every element of ${\cal S}$ is nilpotent and of rank $n-1$, where $J$ is defined in Notation \ref{notaz} and ${\cal Z}$ is a vector subspace.
Let us define $$k= |\{ a \in \{2,\ldots, n-1\}| \; A_{a,1}= 0 \; \forall A \in {\cal Z}\}|$$  
and $${\cal L}= \{A \in {\cal Z} | \; A_{(i)}=0 \; \mbox{\rm for} \; i=2, \ldots, n \}.$$
If for every  $n' \in \{1, \ldots, n-1\}$ 
the maximal dimension of an affine subspace of nilpotent  $n' \times n'$ matrices over $K$ with rank equal to $n'-1$ 
 is less than or equal to 
$ \frac{(n'-1)(n'-2)}{2},$ then  at least one of the following statements holds:

\begin{itemize}
\item $\dim({\cal Z}) \leq  \frac{(n-1)(n-2)}{2}$
\item $ \dim ({\cal L}) \geq k+1$.
\end{itemize}
\end{lemma}

\begin{proof}
We  claim that, if $k=n-2$, then  $A^{(1)}=0$  for any $ A \in {\cal Z}$; in fact, $A_{n,1}=0$ since $A\in {\cal Z}$ (see Lemma \ref{primorem}), and finally $A_{1,1}=0$ holds as well, since $A_{a,1}=0$ for every $  a \in \{2,\, \ldots ,\, n\}$ and $A$ is nilpotent;
so, if we apply  Lemma \ref{triangblocchi} with $r=1$, we can
conclude that $\dim({\cal Z}) \leq \frac{(n-1)(n-2)}{2}$, which completes the argument when $k=n-2$. 

Now let us consider  the case $0\leq k\leq n-3$; 
if $k >0$, let $i_1, \ldots, i_k \in \{2, \ldots,n-1\}$  with $i_1 <\ldots < i_k$ be 
such that $A_{i_1,1}= \ldots= A_{i_k,1}=0 $ for every $A \in {\cal Z}$;
we define the map
$$
\begin{aligned}
    f : {\cal Z} \longrightarrow & K^{n-2-k} \\
    C \longmapsto & 
    (C_{2,1}, \ldots , \widehat{C_{i_1,1}} , \ldots,
    \widehat{C_{i_k,1}} , \ldots , C_{n-1,1}) \, ,
\end{aligned}
$$
where  $ \widehat{C_{i_h,1}} $ means that 
we are omitting $ C_{i_h,1} $ and
 it is understood that in the particular case when $k=0$ the map $f$ sends $C\in {\cal Z}$ to $(C_{2,1}, \ldots , C_{n-1,1})\in K^{n-2}$. 

We notice that, if $C\in \mbox{Ker}(f)$, then $C^{(1)}=0$: in fact, $C\in \mbox{Ker}(f)$ implies that $C_{a,1}=0$ for $a=2,\ldots, n-1$; moreover, $C_{n,1}=0$ holds since $C\in {\cal Z}$ (see Lemma \ref{primorem}); finally the nilpotency of $C$ and the fact that  $C_{a,1}=0$ for $a=2,\ldots, n$ imply the vanishing of $C_{1,1}$ as well.

So, we define the map
$$
\begin{aligned}
    g : \mbox{Ker}(f) \longrightarrow & M((n-1)\times(n-1), K) \\
    C \longmapsto & 
    C_{(2,\ldots,n)}^{(2,\ldots,n)} \; .
\end{aligned}
$$

We claim that $J_{n-1}+\mbox{Im}(g)$ is an affine subspace of $M((n-1)\times(n-1), K)$ of nilpotent matrices of  rank $n-2$. In fact,  each element of $\mbox{Ker}(f) $ can be regarded as an upper triangular block matrix whose diagonal blocks are of size $1\times 1$ and $(n-1)\times (n-1)$, therefore we can argue similarly to the first part of the proof of Lemma \ref{triangblocchi} and conclude that $J_{n-1}+\mbox{Im}(g)$ is an affine subspace of $M((n-1)\times(n-1), K)$ of nilpotent matrices of rank $n-2$.

We claim that $\mbox{Ker}(g) = {\cal L}$. In fact,  $\mbox{Ker}(g)$ consists of matrices $C \in {\cal Z}$ satisfying both  $C^{(1)} = 0$ (since the domain of $g$ is $\mbox{Ker}(f)$) and $C_{(2,\ldots,n)}^{(2,\ldots,n)}=0$; that  is equivalent to $C_{(i)}=0$ for every $ i\in \{2,\ldots n\}$ and $C_{1,1}=0$; so  $\mbox{Ker}(g) \subseteq  {\cal L}$. To prove the other inclusion, observe that, if $C\in {\cal Z}$ satisfies $C_{(i)}=0$ for every $ i\in \{2,\ldots n\}$, then $C_{1,1}=0$ by the nilpotency of $C$. 
Thus we achieved $\mbox{Ker}(g) = {\cal L}$.

Now, in order to complete the argument, let us prove that in case $ \dim ({\cal L}) \leq k$,  the dimension of ${\cal Z}$ is less than or equal to $\frac{(n-1)(n-2)}{2}$.

Suppose  $ \dim ({\cal L}) \leq k$; we have 
$$
\begin{aligned}
    \dim({\cal Z}) &= \dim\mbox{Im}(f) + \dim\mbox{Ker}(f) 
    \leq \dim\mbox{Ker}(f) + n-2-k =\\
    &=\dim\mbox{Ker}(g) + \dim\mbox{Im}(g) + n-2-k \leq \\
    &\leq \dim\mbox{Ker}(g) + \frac{(n-2)(n-3)}{2} + n-2-k \leq \\
    &\leq k + \frac{(n-2)(n-3)}{2} + n-2-k = 
    \frac{(n-1)(n-2)}{2},
\end{aligned}
$$
where in the second-to-last inequality we applied to $\mbox{Im}(g)$  the hypothesis in our statement with $n'=n-1$, and in the last inequality we used the fact that $\mbox{Ker}(g) = {\cal L}$.
\end{proof}

Now we are ready to prove  Theorem \ref{dimker=1}.

\begin{proof}[Proof of Theorem \ref{dimker=1}]
Obviously the affine subspace 
$$\{A \in M(n \times n, K) \; | \; A_{i,j}= 1 \; \mbox{\rm if }  j=i+1 \; \mbox{\rm and } 
A_{i,j}= 0 \; \mbox{\rm if }  j<i+1 \}$$
has dimension $\frac{(n-1)(n-2)}{2}$ and every of its elements is a nilpotent matrix of rank $n-1$.

To prove that the maximal dimension of an affine subspace of nilpotent  $n \times n$ matrices over $K$ with rank equal to $n-1$  is less than or equal to 
$ \frac{(n-1)(n-2)}{2}$ we argue by induction on $n$. The cases $n=1,2$ are easy and left to the reader. So let $n \geq 3$.

Let ${\cal S}$ be  an affine 
subspace of nilpotent  $n \times n$ matrices over $K$ with rank equal to $n-1$. We want to prove that ${\cal S}$ has dimension less than or equal to $\frac{(n-1)(n-2)}{2}$. 

By operating a similarity transformation, we can suppose that $${\cal S}= J +{\cal Z} $$ where ${\cal Z}$ is a vector subspace of $M(n \times n,K)$.

If the first column of every element of ${\cal Z}$ is zero, we can conclude by Lemma \ref{triangblocchi}.
So we can suppose there exists an element of ${\cal Z}$ whose first column is nonzero.

By contradiction suppose that $\dim ({\cal Z}) > \frac{(n-1)(n-2)}{2}$.

Let us define $$k:= |\{i \in \{2, \ldots, n-1\}\; | \; A_{i,1}=0 \; \forall A \in {\cal Z}\}|$$
and 
$$ \overline{l}= \max \{l \in \{1, \ldots , n\}| \;\exists A \in {\cal Z} \; | \;  A_{l,1} \neq 0\}$$
(the set whose maximum we are considering is nonempty for our assumption).

Obviously $k \geq 0$
and $ \overline{l} \geq n-k-1$
(if we had $\overline{l} \leq n-k-2 $, the last $k+2$ elements of the first column of every element of ${\cal Z}$ would be zero, contradicting the definition of $k$). 

Let $\overline{A} \in {\cal Z} $ be such that 
$\overline{A}_{\overline{l},1} \neq 0$.

By Lemma \ref{dimL} we have: 
\begin{equation} \label{dim L} \dim( {\cal L}  ) \geq k+1,
\end{equation} 
where $${\cal L} := \{A \in {\cal Z}  | \; A_{(i)}=0 \; \mbox{\rm for} \; i=2, \ldots, n \}.$$
Let us define 
\begin{equation} \label{overlinei}
\overline{i}= \min \{i \in \{1, \ldots, n\}\; | \; \exists X\in {\cal L}  \; \mbox{\rm s.t.} \; X_{1,i} \neq 0 \};
\end{equation}
observe that the set  whose minimum we are considering is nonempty since $\dim( {\cal L} ) \geq k+1 >0$.
Moreover, since  $\dim( {\cal L} ) \geq k+1$, we must have $\overline{i} \leq n-k$.
Let $\overline{X} \in {\cal L} $ be such that $\overline{X}_{1, \overline{i}} \neq 0$.

\underline{CASE $\overline{l} - \overline{i} \geq 0 $.}

Observe that, by Lemma \ref{primorem}, the elements of   $\langle {\cal S}  \rangle$ (see Notation \ref{notaz}) are nilpotent matrices.
By Proposition \ref{MMS} applied with   ${\cal V}= \langle  {\cal S}  \rangle $, we have that
$$ tr(  (J+s \overline{A})^{m}  \overline{X}) =0 $$
for every positive natural number $m < |K|$ and for every $s \in K$, hence 
$$ \sum_{i= \overline{i}}^{n}
\overline{X}_{1,i}  \left((J+s \overline{A})^{m}\right)_{i,1}=0 
$$ for every positive natural number $m <|K|$ and
for every $s \in K$. The left-hand side of the equality above is obviously a 
polynomial in $s$ of degree less than or equal to $m$; so, if it is nonzero, it has at most $m$ roots.
Hence, if $|K| >m$, all the coefficients of the polynomial must be equal to zero. In particular,
if $|K| >m$,
the coefficient of the term of degree $1$, which is equal to
 $$\sum_{i= \overline{i}}^{n}
\overline{X}_{1,i} \left(\sum_{x,y \in \N \; | \; x+y= m-1 } J^x\overline{A} J^y \right)_{i,1},
$$ 
must be zero;
if $y >0$, we have that $\left(J^x\overline{A} J^y \right)_{i,1} =0$ (see part (i) of Lemma \ref{J}); hence we get that 
the coefficient of the term of degree $1$ is 
$$\sum_{i= \overline{i}}^{ n}
\overline{X}_{1,i} ( J^{m-1}\overline{A})_{i,1},
$$ 
that is 
\begin{equation} \label{xx}
\sum_{i= \overline{i}}^{n-m+1}
\overline{X}_{1,i}  \overline{A}_{i+m-1,1}.
\end{equation}
So we have deduced that,  if $|K| >m$, the sum in $(\ref{xx})$ is zero.
Now take $m=\overline{l} -\overline{i} +1$; observe that it is less than or equal to $n$, which is less than $|K|$ by the assumptions of the theorem.
So, for such choice of $m$,  the sum in (\ref{xx}) must be zero, but, for such choice of $m$,  the sum in (\ref{xx}) 
becomes
 $$
\overline{X}_{1,\overline{i}}  \overline{A}_{\overline{l},1}
$$
since $\overline{A}_{i+m-1,1} = \overline{A}_{i+\overline{l} - \overline{i},1}=0$ for $i \geq \overline{i}+1$.
So we get a contradiction because $
\overline{X}_{1,\overline{i}}  \overline{A}_{\overline{l},1}
$ is nonzero.

\underline{CASE $\overline{l} - \overline{i} < 0 $.}

Since $ \overline{l} \geq n-k-1$ and $\overline{i} \leq n-k$, we get 
$\overline{l} - \overline{i} \geq  n-k-1 -n+k=-1 $, hence we must have $\overline{l} - \overline{i} =-1 $ and  
\begin{equation} \label{lisegnato}
\overline{l} = n-k-1,
\hspace*{2cm}\\
\overline{i} = n-k;
\end{equation}
therefore, for the definition of $\overline{i}$ (see (\ref{overlinei})), we have that $\dim ( {\cal L} ) \leq k+1$ and 
then, by (\ref{dim L}),  $\dim ( {\cal L} ) =k+1$.
So we can deduce that $${\cal L} =\langle E_{1,n-k}, \ldots, E_{1,n}\rangle. $$

Now we want to prove that 
\begin{equation} \label{scopo}
 A_{a,b}=0  \;\;\;\; \mbox{\rm for} \;\; 
 a=n-k, \ldots, n , \;\;\;\; b =1, \ldots, n-k-1, \;\;\;\;\forall A \in {\cal Z}  ;
 \end{equation}
 observe that we already know the claim for $b=1$ since $\overline{l}=n-k-1$.
By contradiction, suppose 
$$ \exists (a,b)  \in \{n-k, \ldots, n\}  \times  \{2, \ldots, n-k-1\} \; | \; \exists A \in {\cal Z}   \; \mbox{\rm with } \; A_{a,b} \neq 0.$$ 

Among the couples $(a,b) $ of this kind, we choose one such that $a-b $ is maximal; we call it 
$(\overline{a}, \overline{b})$; 
 so 
\begin{equation} \label{elen}
X_{a,b}=0 \; \forall (a,b) \in \{ \overline{a}, \ldots, n\} \times  \{1, \ldots, \overline{b}\} -
\{( \overline{a} ,\overline{b})\}, \; \forall X \in {\cal Z} .
\end{equation}
Let $A' \in {\cal Z}   $ be such that 
$A'_{\overline{a}, \overline{b}} \neq 0$.

By our assumption on the cardinality of the field, we can choose $h \in K$ such that $A' + h \overline{A}$ has the entries $(\overline{l},1) $ and $( \overline{a} ,\overline{b})$ nonzero;  for such $h$ we define $$ \hat{A}= A' + h \overline{A}$$

From the definition of $\overline{l}$ and  from (\ref{elen}), we have:  
\begin{equation} \label{el}
 \left\{ \begin{array}{l}
\hat{A}_{i,1}=0 \; \forall i > \overline{l}, \\
\hat{A}_{a,b}=0 \;  \forall (a,b) \in \{{\overline{a}}, \ldots, n\} \times  \{1, \ldots, \overline{b}\} -\{( \overline{a} ,\overline{b})\}. 
\end{array} \right.
\end{equation}

Now observe that, by part (iv) of Lemma \ref{J}, we can choose $c_i \in K$ for $i=1,\ldots, n-1$ such that, if we define
$$ C = I+\sum_{i=1}^{ n-1}  c_i J^i  $$
and
$$ \tilde{A}=C \hat{A} C^{-1} , $$
we have that $\tilde{A}_{i,1}=0$ for every 
$i < \overline{l}$. Furthermore, by part (iii) of Lemma \ref{J}, from (\ref{el}) and from our  choice of $h$, we have:

\begin{equation} \label{elena}
 \left\{ \begin{array}{l}
\tilde{A}_{ \overline{l},1}= \hat{A}_{ \overline{l},1} \neq 0 , \\
\tilde{A}_{  \overline{a} ,\overline{b}}= \hat{A}_{  \overline{a} ,\overline{b}} \neq 0 ,\\
\tilde{A}_{i,1}=0 \; \forall i > \overline{l}, \\
\tilde{A}_{a,b}=0 \; 
\forall (a,b) \in \{{\overline{a}}, \ldots, n\} \times  \{1, \ldots, \overline{b}\} -\{( \overline{a} ,\overline{b})\}. 
\end{array} \right.
\end{equation}

So 
\begin{equation} \label{primacolAtilde}
\tilde{A}^{(1)}= \tilde{A}_{\overline{l},1} e_{\overline{l}}.
\end{equation}

Let $$\tilde{{\cal S} }= C {\cal S} C^{-1}, \;\;\; \tilde{{\cal Z} }= C {\cal Z}  C^{-1}, \;\; \; \tilde{{\cal L} }= C {\cal L}  C^{-1}.$$ 

Obviously $ \tilde{{\cal L} }$ has the same dimension as $ {\cal L} $ and it is contained in ${\cal L} $ (see Lemma \ref{J}, part (i)); so $ \tilde{ {\cal L} }={\cal L} $. Furthermore,  by part (iii) of Lemma \ref{J}, $$\max \{l \in \{1, \ldots , n\}| \;\exists A \in \tilde{{\cal Z} } \; | \;  A_{l,1} \neq 0\}$$ is equal to $\overline{l}$.

Moreover, the space $\tilde{{\cal S} }$ is an affine space of nilpotent matrices of rank $n-1$ and it is equal to $J + \tilde{{\cal Z} }$ since $CJ C^{-1}=J$; by Lemma \ref{primorem}, every element of $\tilde{{\cal Z} }$ is nilpotent, so $ \langle \tilde{ {\cal S} }\rangle$ is  a linear subspace of nilpotent matrices.

Observe that $k < n-2 $, that is $n-k-2 \geq 1$, since we have supposed that there exists an element of ${\cal Z}$ with nonzero first column.
By Proposition \ref{MMS} applied with  $ {\cal V}= \langle \tilde{{\cal S} }\rangle $, we have that, 
 for $r=1, \ldots, k+1$, $t=1, \ldots, n-k-2$,
 $$ tr ((J+s\tilde{A})^{n-k-t} E_{1, n-k-1+r})=0$$ 
for every $s \in K$ (recall that $\tilde{ {\cal L} }={\cal L} = \langle E_{1,n-k}, \ldots, E_{1,n} \rangle $ and that $\tilde{ {\cal L} } \subseteq \tilde{{\cal Z} } \subseteq \langle \tilde{{\cal S} }\rangle $), that is 
 $$ ((J+s\tilde{A})^{n-k-t})_{n-k-1+r,1}=0$$ 
for every $s \in K$.
Obviously  $ ((J+s\tilde{A})^{n-k-t})_{n-k-1+r,1}$ is a polynomial in $s$ of degree less than or equal to $n-k-t $, which is less than or equal to $ n-1$. So,
by the assumption on the cardinality of $K$,
every coefficient of this polynomial must be zero. 
Let us call $\star_{r,t}$ the coefficient of the  term of degree $2$.

We have:
$$ \star_{r,t}=  \sum_{x,y,z \in \N, x+y+z =n-k-t-2}  \left(
J^x \tilde{A} J^y \tilde{A} J^z
\right)_{n-k-1+r,1}= $$
$$= \sum_{x=0}^{ n-k-t-2}  \left(
J^x \tilde{A} J^{n-k-t-2-x} \tilde{A} 
\right)_{n-k-1+r,1},$$ since $\left(
J^x \tilde{A} J^y \tilde{A} J^z
\right)_{n-k-1+r,1} =0 $ if $z >0$.

Hence 
\begin{equation} \label{star}
\begin{split} 
 \star_{r,t} 
 = \sum_{x=0}^{ n-k-t-2}  \left(
J^x \tilde{A}\right)_{(n-k-1+r)} \left( J^{n-k-t-2-x} \tilde{A} 
\right)^{(1)}= \\ 
= \sum_{x=0}^{\min\{ n-k-t-2, k+1-r \}}  
 \tilde{A}_{(n-k-1+r+x)}  \,\tilde{A}_{n-k-1,1} \; e_{1+t+x}= \\ 
= \sum_{x=0}^{\min\{ n-k-t-2, k+1-r \}}    
 \tilde{A}_{n-k-1+r+x,1+t+x} \; \tilde{A}_{n-k-1,1} ,
 \end{split}
 \end{equation}
 where in the second equality we used formula (\ref{primacolAtilde}), the fact that $\overline{l}=n-k-1$ (see (\ref{lisegnato})) and the fact that $ \left(
J^x \tilde{A}\right)_{(n-k-1+r)} =0$ if $x > k+1-r$.
As we have already said, $\star_{r,t}$ must be zero for $r=1, \ldots, k+1$ and $ t =1, \ldots, n-k-2$.

Let us define $$ m_{r,t} =\min\{ n-k-t-2, k+1-r \} ;$$
obviously $$ m_{r+1,t+1} = m_{r,t}-1.$$
Observe that, for  $r=1, \ldots, k$ and $ t =1, \ldots, n-k-3$, we have:
$$ 0= \star_{r,t} - \star_{r+1,t+1} = $$
$$ = \tilde{A}_{n-k-1,1} \left( 
\sum_{x=0}^{ m_{r,t} }  
 \tilde{A}_{n-k-1+r+x,1+t+x}-
\sum_{x=0}^{m_{r+1,t+1} }      
 \tilde{A}_{n-k-1+r+1+x,1+t+1+x}\right)=$$
 $$ = \tilde{A}_{n-k-1,1} \left( 
\sum_{x=0}^{ m_{r,t} }  
 \tilde{A}_{n-k-1+r+x,1+t+x}-
\sum_{x=0}^{m_{r,t}-1 }      
 \tilde{A}_{n-k-1+r+1+x,1+t+1+x}\right)=$$
 $$ = \tilde{A}_{n-k-1,1} \left( 
\sum_{x=0}^{ m_{r,t} }  
 \tilde{A}_{n-k-1+r+x,1+t+x}-
\sum_{x'=1}^{ m_{r,t}
 }   
 \tilde{A}_{n-k-1+r+x',1+t+x'}\right)=$$
$$=\tilde{A}_{n-k-1,1} \; \tilde{A}_{n-k-1+r,1+t}.$$
Since $\tilde{A}_{n-k-1,1} $, that is $\tilde{A}_{\overline{l},1} $ (which is equal to $\hat{A}_{\overline{l},1} $) is nonzero, we get 
\begin{equation} \label{bo1} \tilde{A}_{n-k-1+r,1+t}=0
\end{equation}
 for $r=1, \ldots, k, $ and $ t =1, \ldots, n-k-3$.

Now let us prove the equality (\ref{bo1}) also in the cases $r=k+1$ or $t=n-k-2$.

First let us consider the equality $\star_{r,t}=0$ with $r=k+1$: by (\ref{star})
$$0= \star_{k+1,t}
= \tilde{A}_{n,1+t} \; \tilde{A}_{n-k-1,1} .$$
Since $\tilde{A}_{n-k-1,1} $ is nonzero, we get the equality $(\ref{bo1})$ also for $r=k+1$.
Now  consider the equality $\star_{r,t}=0$ with $t=n-k-2$: by (\ref{star})
$$0= \star_{r,n-k+2}
=   \tilde{A}_{n-k-1+r,n-k-1} \; \tilde{A}_{n-k-1,1} .$$
Since $\tilde{A}_{n-k-1,1} $ is nonzero, we get the equality $(\ref{bo1})$ also for $t=n-k-2$.
Hence we have proved: 
$$ \tilde{A}_{n-k-1+r,1+t}=0 \;\;\;\; \mbox{\rm for} \;\;  r=1, \ldots, k+1, \; \;  \;\;  t =1, \ldots, n-k-2.$$ Thus
$$ \tilde{A}_{a,b}=0 \;\;\;\; \mbox{\rm for} \;\; 
 a=n-k, \ldots, n, \;\;\;\;  b=2, \ldots, n-k-1 , $$
 which contradicts the second statement in (\ref{elena}).
 Hence we proved (\ref{scopo})
and we can conclude
 by Lemma \ref{triangblocchi}.
\end{proof}

We  now consider affine subspaces of nilpotent matrices of rank $1$. Before stating the result, we recall that, if ${\cal V} $ is  a linear space of matrices with rank at most $1$, then all the nonzero matrices in ${\cal V} $ either have a common $1$-dimensional image or have a common $(n-1)$-dimensional kernel, see for instance Lemma 2 in \cite{A-L}.

\begin{proposition} \label{rank=1}
Let $n \in \N$ with $n \geq 2$ and let $K$ be  a field such that $|K|\geq 3$; then the maximal dimension of an affine subspace of nilpotent  $n \times n$ matrices over $K$ with rank equal to $1$  is $n-2$.
\end{proposition}

\begin{proof}
To prove the statement,
 first notice that 
$E_{1,2} + \langle E_{1,3}, \ldots, E_{1,n} \rangle$
is an affine subspace of nilpotent  $n \times n$ matrices with rank equal to $1$  and dimension $n-2$, hence
 the maximal dimension of an affine subspace of nilpotent  $n \times n$ matrices over $K$ with rank equal to $1$  is at least $n-2$.
 
Now, let $ {\cal S} $ be an affine subspace of nilpotent  $n \times n$ matrices over $K$ with rank equal to $1$ and let ${\cal S}  =P+ {\cal Z} $, where ${\cal Z} $ is a linear subspace of $M(n \times n, K) $ and $P \in M(n \times n, K) $. Since $|K| \geq 3$, also the nonzero elements of ${\cal Z} $ have rank $1$ (in fact, for every $A \in {\cal Z} $,  the determinant of every $2 \times 2 $ submatrix of $P +t A$ is  polynomial in $t$ of degree at most $2$ and is zero for every $t \in K$; hence its leading coefficient must be zero). 

Therefore all the nonzero elements of the span ${\cal R} :=\langle \{P\}  \cup {\cal Z}  \rangle $ have rank $1$, hence  
 either  have a common $1$-dimensional image or have a common $(n-1)$-dimensional kernel. 
 
 Moreover, observe that a nilpotent matrix of rank $1$ has square zero, so ${\cal S}$ is an affine space of square zero matrices. Hence $(P+tA)^2=0$ for every $A \in {\cal Z}$ and for every $t \in K$ and, arguing as in the proof of Lemma \ref{primorem}, we get that $A^2=0$  for every $A \in {\cal Z}$.
 Consequently all the elements of ${\cal R}$ have square zero; in particular,   all the elements of ${\cal R}$ are nilpotent.
 
 Since   all the nonzero elements of ${\cal R} $ either  have a common $1$-dimensional image or have a common $(n-1)$-dimensional kernel and all the elements of ${\cal R}$ are nilpotent,
 we have that  $ \dim ({\cal R}) \leq n-1$.

 Finally, observe that  $P \not \in {\cal Z} $ (otherwise the null matrix would be in ${\cal S}$, which is impossible because every element of ${\cal S}$ has rank $1$); hence  we have that $$\dim( {\cal S} ) = \dim({\cal Z} )= \dim({\cal R} )-1 \leq n-2,$$ as we wanted to prove. 
\end{proof}

\begin{remark}\label{counterexample}
    We provide an example to show that Theorem \ref{dimker=1} does not hold if we drop the assumption that $|K| 
    \geq n+1.$
    Moreover, the same example also shows that Proposition \ref{rank=1} does not hold if we weaken the assumption $|K|\geq 3$.
    
    In fact, we consider $n=2$ and $K=\mathbb{Z}_2$,
    and we observe that 
    $$
    {\cal Z} =
    \left\{
    \left(
\begin{array}{cc}
    0 & a\\
    a & 0 
\end{array} \right)
\; | \; a\in\mathbb{Z}_2
    \right\}
    $$
    despite containing a matrix that is not nilpotent, has the property that
    ${\cal S} :=J+{\cal Z} $ is an affine subspace of nilpotent matrices with rank equal to $1$ and ${\cal S} $ has dimension $1$.  
\end{remark}

We finish with a conjecture on the maximal dimension of affine subspace of nilpotent matrices of constant  rank:

\begin{conjecture} \label{conjecture}
     Let $n \in \N-\{0\}$ and
let $K$ be a field. 
Let ${\cal S} $ be an affine subspace of $M(n \times n, K)$ such that every element of ${\cal S} $ is nilpotent, and assume that each element of ${\cal S} $ has rank equal to $r\geq 1$. Then we conjecture that, if $|K|$ is sufficiently large, then the maximal dimension ${\cal S} $ can attain is
$$
\sum_{i=1}^r (n-1-i) = \frac{r}{2}(2n-r-3)\,.
$$
\end{conjecture}

Observe that 
Theorem \ref{dimker=1} and Proposition \ref{rank=1} confirm the conjecture when $r=n-1$ and $r=1$ respectively.
Observe also that the affine subspace 
$$\{A \in M(n \times n, K) \; | \; A_{i,i+1}= 1 \; \mbox{\rm for }  i=1,\ldots, r \; \mbox{\rm and } 
A_{i,j}= 0 \; \mbox{\rm if }\;  j \leq i \; \mbox{\rm or}  \;i \geq r+1\}$$
has dimension $ \frac{r}{2}(2n-r-3)$ and every of its elements is a nilpotent matrix of rank $r$.

\bigskip

{\bf Acknowledgments.}
In the first version of the paper the assumption on the cardinality of the field was stronger. De Seguins Pazzis 
let us know that Theorem \ref{dimker=1} could be obtained 
also with another strategy using the ideas in \cite{DP5} and the assumption on the cardinality of the field weakened. After his communication, 
we succeeded in weakening the assumption on the field by making some changes and in  simplifying our original proof.  
So we are indebted with De Seguins Pazzis for his precious communication, which stimulated us to improve the paper, and  we thank him warmly.

Finally, 
we  wish to express our sincere gratitude 
 to the anonymous reviewer  for his/her valuable suggestions, which helped us to improve the paper very much; in particular, we thank him/her for the suggestion on how to simplify the proof of Proposition \ref{rank=1}.

\end{document}